\documentclass{amsart}
\usepackage{amsfonts,amssymb,amscd,amsmath,enumerate,verbatim,calc}
\usepackage{xypic}

\newcommand{\CM}{Cohen-Macaulay}

\newcommand{\GCM}{generalized Cohen-Macaulay}

\newcommand{\ff}{\text{if and only if}}
\newcommand{\wrt}{with respect to}
\newcommand{\B}{\mathcal{B} }
\newcommand{\aF}{\mathfrak{a} }

\newcommand{\n}{\mathfrak{n} }
\newcommand{\m}{\mathfrak{m} }
\newcommand{\M}{\mathfrak{M} }
\newcommand{\q}{\mathfrak{q} }
\newcommand{\A}{\mathfrak{a} }
\newcommand{\R}{\mathcal{R} }
\newcommand{\Z}{\mathbb{Z} }

\newcommand{\rt}{\rightarrow}
\newcommand{\xar}{\longrightarrow}
\newcommand{\ov}{\overline}
\newcommand{\sub}{\subseteq}

\newcommand{\wt}{\widetilde }

\newcommand{\image}{\operatorname{image}}

\newcommand{\reg}{\operatorname{reg}}
\newcommand{\grade}{\operatorname{grade}}
\newcommand{\depth}{\operatorname{depth}}
\newcommand{\e}{\operatorname{end}}
\newcommand{\ann}{\operatorname{ann}}

\newcommand{\amp}{\operatorname{amp}}
\newcommand{\Hom}{\operatorname{Hom}}
\newcommand{\Ext}{\operatorname{Ext}}

\theoremstyle{plain}

\newtheorem{theorem}{Theorem}[section]
\newtheorem{corollary}[theorem]{Corollary}
\newtheorem{lemma}[theorem]{Lemma}
\newtheorem{proposition}[theorem]{Proposition}

\theoremstyle{definition}
\newtheorem{definition}[theorem]{Definition}
\newtheorem{s}[theorem]{}
\newtheorem{remark}[theorem]{Remark}
\newtheorem{example}[theorem]{Example}
\newtheorem{observation}[theorem]{Observation}

\theoremstyle{remark}

\begin{document}

\title[Higher Associated graded modules]{Ratliff-Rush Filtration, regularity
 and
 \\ depth of Higher Associated graded modules \\ Part I}
\author{Tony~J.~Puthenpurakal}
\date{\today}
\address{Department of Mathematics, IIT Bombay, Powai, Mumbai 400 076}

\email{tputhen@math.iitb.ac.in}
 \begin{abstract}
In this paper we introduce a new technique to study associated graded modules. Let $(A,\m)$ be a
Noetherian local ring with $\depth A \geq 2$.  Our techniques gives a necessary and sufficient condition
for $\depth G_{\m^n}(A) \geq 2$ for all $n \gg 0$.
Other applications are also included; most notable is an upper bound regarding the Ratliff-Rush
filtration.
\end{abstract}
 \maketitle
\section*{introduction}

Let $(A, \m)$ be a Noetherian local ring of dimension $d$ with
residue field $k = A/\m$. Let $M$ be a finite, that is to say, finitely
 generated
$A$-module of dimension $r$ and let  $I$ be an \emph{ ideal of definition} for $M$ i.e.,
   $\lambda(M/IM)$ is
 finite (here $\lambda(-)$ denotes length).
The \emph{Hilbert function }of $M$ \wrt \ $I$ is $H^{I}(M,n) =  \lambda(I^nM/I^{n+1}M)$ for $n \geq 0$.

When $M = A$ is \CM \ a fruitful area of research has been to study the interplay between Hilbert functions and
properties of the blowup algebra's of $A$ \wrt \ $I$, namely, the \emph{Rees ring } $R(I) = \bigoplus_{n \geq 0}I^nt^n$, the \emph{ extended Rees ring } $R(I)^* = \bigoplus_{n \in \Z}I^nt^n$ (here $I^n = A$ for $n < 0$) and  the
\emph{associated graded ring} $G_I(A) = \bigoplus_{n \geq 0}I^n/I^{n+1}$.
 See the text's
 \cite[Section 6]{VaSix} and  \cite[Chapter 5]{VasBook} for nice surveys on this subject. Graded local cohomolgy  has played an important role in this subject.
 For  various applications
   see  \cite[4.4.3]{BH},\cite{Durham}, \cite{VerJoh},
\cite{Blanc},  \cite{ItN},
   and \cite{HMc}.
In this paper we introduce a \emph{new} technique to study some questions in this area.

\s\label{tech}
\textbf{Technique:}
We study $L^I(M)  = \bigoplus_{n \geq 0} M/I^{n+1}M$; a \emph{not finitely generated}  $R(I)$-module.
Set $\M = \m \oplus R(I)_+$.
It has the following
properties.
\begin{enumerate}[\rm 1.]
\item
For $0 \leq i \leq  \depth M - 1$ the local cohomology modules $H^{i}_{\M}(L^I(M))$ are  *-Artinian.
 Furthermore for each $i = 0, \ldots, \depth M -1$, $\lambda(H^{i}_{\M}(L^I(M))_n)$ is finite for all $n \in \mathbb{Z}$ and it coincides with a polynomial for all $n \ll 0$.
\item
$L^I(M)(-1)$ behaves well \wrt \  the Veronese functor. Clearly
\[
L^I(M)(-1)^{<l>} = L^{I^l}(M)(-1).
\]
\item
Let $\depth M > 0$ and let $x$ be $M$-superficial \wrt \  $I$. Set $N = M/xM$ and $u =xt \in R(I)$. Then
$L^I(M)/uL^I(M) = L^I(N)$.
\end{enumerate}
Before we state the applications we need some notation:

\noindent Let
$G_{I}(M) = \bigoplus_{n \geq 0}I^nM/I^{n+1}M$ be  the \emph{associated
 graded module} of $M$ \wrt \ $I$, considered as a $G_I(A)$ module.
 The ring $G_{I}(A)$ has a  unique graded maximal ideal
$\M_G = \m/I \oplus_{n\geq 1}I^n/I^{n+1} $.
 Set $ \depth G_I(M) =  \grade(\M_G,G_I(M))$.

\noindent \textbf{Applications:}

\noindent \textbf{I.} In \cite[2.2]{EliHA} Elias
proves
 $\depth G_{I^n}(A)$ is constant for $n \gg 0$ (here $I$ \emph{need not be}
$\m$-primary). We call this number $\xi_I(A)$.
It's  to see that
 $I$ has a regular element  if and only if
 $\xi_I(A) \geq 1 $.   However there is no general criteria for
$\xi_I(A) \geq 2$.  In this paper we give  necessary and sufficient conditions for
$\xi_{\m}(A) \geq 2$.
 If $\depth G_{\m^n}(A) \geq 2$ for some $n$ then clearly $\depth A \geq 2$.
  However there are
\CM \ local rings $A$ of dimension two with $\xi_{\m}(A) = 1$ (see example \ref{2dimnotgcm}).

\noindent\textbf{Theorem \ref{geq2b}}
\textit{Let $(A,\m)$ be a Noetherian \ local ring with depth $\geq 2$ and residue field $k$. Let
 $G_{\m}(A) = R/\q$, where $R =  k[X_1,\ldots,X_s]$. Then}

$\displaystyle{\Ext^{s-1}_{R}(G_{\m}(A),R)_{-(s-1)} = 0 \quad \ \textit{if and only if} \quad \depth G_{\m^n}(A) \geq 2 \ \textit{for all} \
n \gg 0.}$

\noindent If $s$ is not too big then this
criteria can be checked
with a computer algebra program. Notice we do not make any assumptions
on the residue field.

\noindent \textbf{II.}
We extend  Elias's result to modules in a special case when $\lambda(M/IM)$ is finite
  i.e., we prove $ \depth G_{I^n}(M)$ is constant for all $  n \gg 0  $ (see  \ref{asyym}).
We call this number $\xi_I(M) $.
Our techniques yield a theoretic way  to check  $\xi_I(M)$. We  use  it
  to construct examples with  \CM  \ local ring $A_r$ of dimension $d =  r +s$ with $r \geq 1, s\geq 1$, $I$ an $\m$-primary ideal
 with $\xi_I(A_r) = r$,  see Example \ref{arbitxi}.
Theorem \ref{geq2b} follows from a more general
criteria (see Proposition \ref{xi2firstprop}) to ensure $\xi_I(M) \geq 2$ for all $n \gg 0$.
Unfortunately this criteria   is not  verifiable  with a computer.
 
Our method to show $\xi_I(M)$ is a constant also indicates a method to
attack a different problem. If $(A,\m)$ is Noetherian local and if $K$ is an ideal in $A$ then recall that
\emph{the fibercone} of $K$ is the graded ring $F(K) = \bigoplus_{n\geq 0}K^n/\m K^n$.
In Theorem \ref{fiberasyym} we show $\depth F(K^n)$ is constant for all $n \gg 0$.

\noindent \textbf{III.}
If $G_{I^s}(M)$ is \CM \ for some $s \geq 1$ then  $G_{I^{sm}}(M)$ is
\CM \ for all $m \geq 1$. So we get $\xi_I(M) = r$, i.e., $G_{I^n}(M)$ is \CM \ for all $n \gg 0$.
In Proposition \ref{pgcm} we show

 If $G_{I^n}(M)$ is \CM \ for some $n \geq 1$
then $G_I(M)$ is generalized \CM \ i.e., $ H^{i}(G_I(M))$, the $i'th$ \emph{local cohomology}
modules of $G_I(M)$ \wrt \  the $*$- maximal ideal  $G_I(A)$, has finite length for $i = 0,\ldots, r -1$.

When $(A,\m)$ is \CM \  and $I$ is  a normal $\m$-primary ideal then Huneke and
 Huckaba  \cite[Corollary 3.8]{HHu} show $\xi_I(A) \geq 2$.
  Our techniques yield a simpler proof of this fact. In Theorem \ref{mynor} we show
   $\depth G_{I^n}(A) \geq 2$ for all $n \geq \reg(G_I(A))$

\noindent \textbf{IV.} The motivation to consider $L^I(M)$ was to
understand certain
aspects of
Ratliff-Rush Filtration's
 (see \cite{RR}).  It is useful
to extend this notion to modules. Define
\[
 \wt{IM} = \bigcup_{j \geq 0}(I^{j+1}M \colon_M I^{j}).
\]
If  $\grade(I,M) > 0$ one can easily show  $\wt{I^n M} = I^nM$ for all $n \gg 0$. It is of some interest to find an
upper bound for
$$  \rho^{I}(M) := \min\{ n \mid \wt{I^iM} = I^iM \ \text{for all} \ i \geq n  \}.$$
In the case $I = \m$
this has applications to
 in certain questions in homological dimensions of
$M/\m^nM$ and $\m^nM$ (see \cite{JavPuArxiv},\cite{Pu4}). In Theorem \ref{mthm}
we prove
\[
 \rho^{I}(M) \leq \max\{0, a_{1}\left(G_{I}(M)\right)+1  \}\quad \text{if} \ \grade(I,M)  > 0.
\]

Assume $\depth M \geq 2$ and let $x$ is $M$-superficial \wrt \ $I$. Set $N = M/xM$.
 Clearly $\ov{\wt{I^sM}} \subseteq \wt{I^sN}$ for each $s \geq 1$.
We have the following \emph{natural} exact sequence:
\[
 0 \xar \frac{(I^{n+1}M\colon_M x)}{I^nM} \xar \frac{\wt{I^nM}}{I^nM} \xrightarrow{\alpha_{n-1}^{x}}\frac{\wt{I^{n+1}M}}{I^{n+1}M} \xrightarrow{\rho_n} \frac{\wt{I^{n+1}N}}{I^{n+1}N}.
\]
Here $\alpha_{n-1}^{x}(p + I^nM) = xp + I^{n+1}M$ and $\rho_n$ is the natural map.
 Another motivation for me was to extend
 the  above exact sequence. We do it in \ref{longHpara}.

\noindent \textbf{V.}
Our techniques enable us to generalize some previously known results (our proofs are simpler
too!)
  In \ref{HMR} we extend
to modules a result of Huckaba and Marley \cite{HMr} relating depths of $R(I)$ and $G(I)$.
In \ref{marextn}
 we generalize
to modules a result  due to Marley \cite{HMr} regarding local cohomology
modules of $G_I(A)$. Our generalization is different than by Hoa. See \ref{marextnD}.

Here is an overview of the contents of the paper. In section one we discuss a few
 preliminaries. We also define Ratliff-Rush filtration for modules and discuss a few of its
properties. In section two we discuss local cohomology of graded
modules over standard algebras over a local ring $R_0$,
particulary when $R_0$ is Artinian. We also prove that if $K$ is
an ideal in a Noetherian local ring then $\depth F(K^n)$ is
constant for all $n \gg 0$. In section 3 we introduce $L^I(M)$ and
prove $H^{i}(L^I(M))$ is *-Artinian for $i = 0,\ldots, \depth M
-1$. We also compute $H^0(L^I(M))$. In section 4 we prove the
upper bound on $\rho^I(M)$. The results stated in \textbf{V.} are
also proved in this section. In section 5 we discuss the behavior
of Ratliff-Rush filtration mpdulo a superficial element.
In section 6 we  investigate the cohomological consequences of
\ref{tech}.3. We also prove \ref{tech}.1. In section 7 we
investigate the cohomological consequences of \ref{tech}.2. We
prove all the results stated in \textbf{II.} and \textbf{III.} As
an  application, we prove a curious result regarding Ratliff-Rush
filtration for ideals, see \ref{curiousm}, \ref{curiousnormal}. In
section 8 we make a critical observation which is used later in
section 9 and in part II of our paper. In section 9 we give
necessary and sufficient criteria for $\depth G_{\m^n}(A) \geq 2$
for all $n \gg 0$.

\section{Preliminaries}
In this paper all rings are commutative Noetherian and all modules
are assumed finite. We  use  terminology from  \cite{BH}.  Let
$(A,\m)$ be a local ring of dimension $d$ with residue field $k =
A/\m$. Let $M$ be an $A$-module. Let $I$ be an ideal in $A$ (not
necessarily an ideal of definition for $M$).
If $p \in M$ is non-zero and  $j$ is the largest integer such that $p \in I^{j}M$,
then we let $p^*$ denote the image of $p$
 in $I^{j}M/I^{j+1}M$. Set $0^* = 0$.

\begin{remark}
Let $x_1,... ,x_s$ be a sequence in $A$ with $x_i \in I$ and
set $J= (x_1,... ,x_s)$.
Set $B = A/J$,  $K = I/J$ and  $N = M/JM$. Notice
\[
G_I(N) = G_K(N) \quad \text{and} \quad \depth_{ G_I(A)} G_I(N) = \depth_{G_K(B)} G_K(N).
\]
\end{remark}

\s For definition and few basic properties of superficial sequences see \cite{Pu1}, pages 86-87.

\s \textbf{Base change:}
\label{AtoA'}
 Let $\phi \colon (A,\m) \rt (A',\m')$ be a local ring homomorphism. Assume
 either $A'$ is a quotient of $A$ or  $A'$ is a faithfully flat $A$
algebra with $\m A' = \m'$. Set $I' = IA'$ and if
 $N$ is an $A$-module set $N' = N\otimes A'$.
 In these case's it can be seen that

\begin{enumerate}[\rm (1)]
\item
$\lambda(N) = \lambda(N')$.
\item
 $H^I(M,n) = H^{I'}(M',n)$ for all $n \geq 0$.
\item
$\dim M = \dim M'$ and  $\grade(K,M) = \grade(KA',M')$ for any ideal $K$ of $A$.
\item
$\depth G_{I}(M) = \depth G_{I'}(M')$.
\end{enumerate}

 \noindent The specific base changes we do are the following:

(i) If $I$ is an ideal of definition of $M$ but $I$ is \emph{not} $\m$-primary then we set $A' =
A/\ann(M)$. Then $I' = (I+ \ann(M))/ \ann(M)$ is $\m'$-primary.
Furthermore $M' = M$ as $A$-modules.

(ii) $A' = A[X_1,\ldots,X_n]_S$ where $S =  A[X_1,\ldots,X_n]\setminus \m A[X_1,\ldots,X_n]$.
The maximal ideal of $A'$ is $\n = \m A'$.
The residue
field of $A'$ is $l = k(X_1,\ldots,X_n)$. Notice that if $I$ is integrally closed then $I'$ is 
also integrally closed.

(a) If the residue field is finite we make this base change with $n =1$ just to ensure existence of
superficial elements.

(b) When $\dim A \geq 2$,
Itoh \cite[Lemma 11]{It} shows that
 there exists a superficial element $y \in I'$ such that the
$A'/(y)$ ideal $I'/(y)$
is  integrally closed ideal.

 (c) When $\dim A \geq 2$ and $I$ is normal,
Itoh   \cite[Theorem 1]{ItN} shows that
 there exists a superficial element $y \in I'$ such that
the $A'/(y)$ ideal $J = I'/(y)$ is asymptotically normal i.e.,
$J^n$ is integrally closed for all $n \gg 0$.

\noindent\textbf{Ratliff-Rush Filtration for Modules:} The notion of Ratliff-Rush filtration of an ideal has proved to be an important technique in the
study of blowup algebras. We extend it to modules.

\begin{definition}
  Consider the following chain of submodules of
$M$:
\[
IM \sub (I^2M\colon_M I) \sub (I^3M\colon_M I^2) \sub \ldots \sub(I^{n+1}M \colon_M I^n)\sub \ldots
\]
Since $M$ is Noetherian this chain of submodules stabilizes. We
denote the stable value to be $\widetilde{IM}$.
We call $\wt{IM}$ to be the \emph{Ratliff-Rush submodule of $M$ associated with $I$}.
 The filtration
$\{\wt{I^nM}\}_{n \geq 1}$ is called the \textit{Ratliff-Rush filtration} of $M$
\wrt \ $I$.
\end{definition}

The next theorem collects the two
most important properties of Ratliff-Rush filtration's. The proof in \cite{RR} in the case of rings extends
to modules.
\begin{theorem}
\label{rrf}
Let $A$ be a ring,
  $M$ an $A$-module and
 $I$  an ideal of $A$. If

\noindent $\grade(I,M) > 0$ then the following holds
\begin{enumerate}[\rm 1.]
\item
$\wt{I^nM} = I^nM$ for all $n \gg 0$.
\item
If $x \in I$ is  $M$-superficial \wrt \ $I$ then
$(\wt{I^{n+1}M}\colon_M x) = \wt{I^nM}$ for all $n \geq 1$. \qed
\end{enumerate}
\end{theorem}

\noindent\textbf{Ratliff-Rush Filtration and base change:}
\begin{observation}
\label{rrhom}
Let $B = A/\q$ for some ideal $\q$ and let $N$ be a finite $B$-module. If
$K = (I + \q)/\q$ is an ideal in $B$, with $I$ an ideal in $A$, then
$\wt{KN} = \wt{IN}$.
In particular  $\wt{K^jN} = \wt{I^jN}$
for each $j \geq 1$.
\end{observation}

Next we deal with flat base change.
\begin{proposition}
\label{rrflat}
Let $(A,\m) \rt (B,\n)$ be a flat homomorphism of local rings. Let $M$ be an $A$-module and let $I$ be an ideal in $A$. Set $M' = M\otimes_AB$ and let $J = IB$.
We then have:
\begin{enumerate}[\rm 1.]
\item
$I^{n+1}M \otimes_A B \cong J^{n+1}M'$ for all $n \geq 0$.
\item
$(I^{n+1}M\colon_M I^n)\otimes B = (J^{n+1}M'\colon_{M'} J^n)$.
\item
$\wt{IM}\otimes_A B = \wt{JM}$.
\item
$\wt{I^nM}\otimes_A B = \wt{J^nM}$ for all $n \geq 1$.
\end{enumerate}
\end{proposition}
\begin{proof}
Since $B$ is  a flat $A$-module   we have $J \cong I\otimes B$.
The statement  $1.$ is a  well-known fact.
For $2.$ see \cite[Theorem 7.4]{Ma}.
 From $2.$  we easily get
3. For 4. fix $n \geq 1$. Since $I^nB = J^n$ we get $\wt{I^nM}\otimes_A B = \wt{J^nM}$ by  $3.$
 \end{proof}

\s
If  $\grade(I,M) > 0$ then  $\wt{I^n M} = I^nM$ for all $n \gg 0$.
This motivates the following definition:
\[
  \rho^{I}(M) := \min\{ n \mid \wt{I^iM} = I^iM \ \text{for all} \ i \geq n  \}.
\]

\section{Polynomial growth of graded local cohomology }
As a reference for local cohomolgy we use  \cite{BSh}.
Let $(R_0,\m_0)$ be a local ring. We say
  $R = \bigoplus_{i \geq 0}R_i$ is a \emph{standard graded} $R_0$-algebra if    $R$ is generated over $R_0$ by finitely
many elements of degree $1$.  Set $R_{+} = \bigoplus_{i \geq
1}R_i$ and $\M_R = \m_0\oplus R_+$. We drop the subscript $R$ if
the ring is clear from the context.
Let $L$ be a (not necessarily finitely generated) graded $R$-module. Define
  $\e(L) = \sup\{n \in \Z \mid L_n \neq 0\}$.
If $\A$ is a  homogeneous ideal in $R$
then we set $H^{i}_{\A}(L)$ to be the \emph{$i$-th local cohomology} module of $L$
\wrt \ $\A$.

\s If $E$ is a finite $R$-module then for each $i \geq 0$ we have $H^{i}_{R_+}(E)_n = 0$ for all $n \gg 0$
(cf. \cite[15.1.5]{BSh}). Set
\begin{align*}
a_{i}(E) &= \e \left( H^{i}_{R_+}(E)\right), \\
\reg E &= \max \{ a_i(E) + i \mid 0\leq i \leq \dim E \}.
\end{align*}
The number $\reg E$ is called the
\emph{(Castelnuovo-Mumford) regularity} of $E$.

\s \label{2artin}  A graded $R$ module $L$ is said to be $*$-\emph{Artinian} if
every descending chain of graded submodules of $L$ terminates.
If $L$ is *-Artinian then  $L_n = 0$ for all $n \gg 0$.
If $M$ is a finite $R$-module then  $H^{i}_{\M}(M)$ is $*$-Artinian for each
$i \geq 0$ \cite[3.6.19]{BH}.

\s\label{verEs}
 Let $l$ be a positive integer. Let $R^{<l>}= \bigoplus_{n\geq 0}R_{nl}$ be the
\emph{$l$-th Veronesean subring} of $R$.
Notice $R^{<l>}$ is also a standard graded $R_0$-algebra and
$\M^{<l>}$ is  the unique graded maximal ideal of $R^{<l>}$.
 If $E$ is a graded $R$-module
and $l$ is a positive integer then the $l'th$ Veronesean submodule of $E$ is
$E^{<l>} :=  \bigoplus_{n \in \Z}E_{nl}$. Clearly $E^{<l>}$ is a graded $R^{<l>}$ module.
The Veronese functor  commutes with local
cohomology cf. \cite[Proposition 2.5 ]{HMc}; i.e. if $\A$ is a homogeneous ideal in $R$ then
\begin{equation}
\label{verE}
H^{i}_{\A^{<l>}}\left(E^{<l>} \right) =  \left(H^{i}_{\A}(E) \right)^{<l>} \quad \text{ for all $i \geq 0$.}
\end{equation}

\s \textbf{Pertinent Examples:}\label{perexam} Let $A$ be local and let $I$ be an ideal.
The Rees ring $R(I) = A[It]$ is a standard graded $A$-algebra. The associated graded
ring is a standard graded $A/I$-algebra. Also $F(I) =  \bigoplus_{n\geq 0}I^n/\m I^n$, the  \emph{fiber cone} of $A$ \wrt \
$I$ is the standard graded $k$-algebra.

The next Lemma is known but not so well-known.

\begin{lemma}
\label{polyartin}
Let $(R_0,\m_0)$ be an Artinian local ring and let $R = \bigoplus_{n \geq 0}R_n$ be a
standard graded $R_0$-algebra. Let $L= \bigoplus_{n \in \mathbb{Z}}L_n$ be an *-Artinian
$R$-module. Then $\lambda(L_n)$ is finite for all $n \in \Z$ and there
exists a polynomial $p_L(t) \in \mathbb{Q}[t]$ with $p_L(n) =\lambda(L_n)$ for
all $n \ll 0$.
\end{lemma}
\begin{proof}
This Lemma follows  from graded Matlis duality cf. \cite[Section 3.6]{BH}. Notice
$R$ is *-complete. Set $k = R_0/\m_0$ and $E = E_{R_0}(k)$ the injective hull of $k$ as
an $R_0$-module. Set
\[
L^{\vee} =\  ^* \Hom_{R_0}(L, E), \quad \text{notice}\  (L^{\vee})_n =  \Hom_{R_0}(L_{-n}, E).
\]
 $L^{\vee}$ is a $R$-module and it is the Matlis dual of $L$.
Since $L$ is *-Artinian,  $L^{\vee}$ is a finitely generated graded $R$-module.
Thus $\lambda((L^{\vee})_n)$ is finite for all $n$.
Notice $\lambda(L_n) = \lambda\left((L^{\vee})_{-n}\right)$ for all $n \in \Z$.
 Let $q(t)$ be the Hilbert polynomial of
$L^{\vee}$. Set $p_L(t) = q(-t)$. Then for $n \ll 0$
\[
p_L(n) = q(-n) = \lambda\left( (L^{\vee})_{-n}\right) = \lambda(L_n).
\]
\end{proof}

\begin{remark}\label{polyartinloc}
\emph{(with the notation as in Lemma \ref{polyartin})}  In particular if $M$ is finite
$R$-module then all $H^{i}_{\M}(M)$ are *-Artinian and so by Lemma  above we get  $\lambda(H^{i}_{\M}(M))_n$ is polynomial
for all $n \ll 0$. For more details see \cite[17.1.9]{BSh}.
\end{remark}

Next we investigate asymptotic depth of Veronese submodules.
\begin{theorem}
\label{asympveron}
Let $(R_0,\m_0)$ be an Artin local ring and let $R = \bigoplus_{n \geq 0}R_n$
be a finitely generated standard $R_0$-algebra. Let $M$ be a finite graded
$R$-module. Then
\[
\depth_{R^{<l>}} M^{<l>} \quad \text{is constant for all} \ n \gg 0.
\]
\end{theorem}
\begin{proof}
Let $\M$ be the $*$-maximal ideal of $R$. Notice that $\M^{<l>}$ is the
$*$-maximal ideal of $R^{<l>}$.
\begin{align*}
\text{Set} \ \  \xi(M) &= \min\{ i \mid H^{i}(M)_{0} \neq 0 \ \text{or} \ \lambda\left(H^{i}(M)\right) = \infty \}, \ \text{and} \\
\amp(M) &= \max \{ \ |n| \ \mid H^{i}(M)_{n} \neq 0 \ \text{for} \ i = 0,\ldots,\xi(M)-1 \}.
\end{align*}
We claim that $\depth M^{<l>} = \xi(M)$ for all $l > \amp(M)$.
Fix $l > \amp(M)$. Notice for $i = 0,\ldots, \xi(M) -1$ we have
\[
H^{i}_{\M^{<l>}}\left(M^{<l>}\right) = H^{i}(M)^{<l>} = 0.
\]
So we have $\depth M^{<l>} \geq \xi(M)$.

\noindent \textbf{Comment:} \emph{Till now our arguments would
work over any non-negatively graded not-necessarily standard algebra over any
local ring. Our final assertion however only works over standard algebra's over Artin local rings.}

Suppose if possible $\depth M^{<l>} > \xi(M)$ for some $l$. Set $s = \xi(M)$. We have
\begin{equation*}
H^{s}(M)_{nl} = H^{s}\left(M^{<l>}\right)_n = 0 \quad \text{for all} \ n \in \Z. \tag{*}
\end{equation*}
So $H^s(M)_0 = 0$.
Since $\lambda(H^{s}(M)_{n})$ is polynomial for all $n << 0$, see  \ref{polyartinloc},
 we have
$H^s(M)_n = 0$ for all $n \ll 0$. Also by \ref{2artin} we have  $H^s(M)_n = 0$ for all $n \gg 0$.
Thus $\lambda(H^s(M))$ is finite.
 So we get $s =\xi(M) \geq s+1$, a contradiction.

Thus
$\depth M^{<l>} = \xi(M)$
for all $l \geq \amp(M)$.
\end{proof}

An immediate application is to asymptotic depth of fibercones.
We have:
\begin{theorem}
\label{fiberasyym}
Let $(A,\m)$ be a Noetherian local ring and let $K$ be an ideal in $A$.
Then $\depth F(K^n)$ is constant for all $n \gg 0$.
\end{theorem}
\begin{proof}
Observe that $F(K^n) = F(K)^{<n>}$. We get the result by Theorem \ref{asympveron}.
\end{proof}

\section{$L^{I}(M)$  }
\s \textbf{Setup:}
\label{setupgen}
In this section $(A, \m)$ is a Noetherian local ring,
$M$ a finite
 $A$-module of dimension $r$.  In this section  the ideal  $I$  \emph{is not
 necessarily an ideal of definition for $M$}. We assume $\grade(I,M) > 0$.

\s
 Set $R(I,M) = \bigoplus_{n \geq 0}I^nM$ the \emph{Rees module} of $M$ \wrt \ $I$.
Clearly $R(I,M)$ is
a finite $R(I)$-module.

\begin{definition}
\label{modstruc}
Set $L^I(M) = \bigoplus_{n \geq 0} M/I^{n+1}M$.
The $A$-module $L^I(M)$ can be given an $R(I)$-module structure as follows.
The Rees ring $R(I)$ is a subring of $A[t]$ and so $A[t]$ is an $R(I)$-module.
Therefore $M[t] = M\otimes_A A[t]$ is an $R(I)$-module. The exact sequence
\[
0 \xar R(I,M) \xar M[t] \xar L^I(M)(-1) \xar 0
\]
defines an $R(I)$-module structure on $ L^I(M)(-1)$ and so on $ L^I(M)$.
\end{definition}

\begin{remark}
\label{mt}
If $x_1,\ldots,x_l \in \m$ is an $M$-sequence then considering $x_1,\ldots,x_l$ in $R_0$, it becomes an $M[t]$ sequence.
So $\Ext^{i}_{R(I)}(R(I)/\M^s,M[t]) = 0$ for all $s \geq 1$ and $0\leq i \leq l-1$.
Therefore $ H^{i}_{\M}\left(M[t]\right) = 0$ for all  $0\leq i \leq l-1$.
\end{remark}

\begin{proposition}
\label{fgenExt}
(Hypothesis as in \ref{setupgen}.)
Set $l = \depth M$, $E =  L^I(M)(-1)$ and $\M = \M_{R(I)}$.
Then
\begin{enumerate}[\rm (a)]
\item
$ H^{i}_{\M}\left(E\right) =  H^{i+1}_{\M}\left(R(I,M)\right)$ for  $0\leq i \leq l-2$.
\item
$ H^{l-1}_{\M}\left(E\right)$ is
 isomorphic to a submodule of $ H^{l}_{\M}\left(R(I,M)\right)$.
\end{enumerate}
In particular $ H^{i}_{\M}\left( L^I(M)\right)$ is *-Artinian for $0 \leq i \leq l -1$.
\end{proposition}
\begin{proof}
When $l = 0$ we have nothing to prove. So assume $ l = \depth M \geq 1$.
 Using Remark \ref{mt}
we get  $ H^{i}_{\M}\left(M[t]\right) = 0$ for all  $0\leq i \leq l-1$.
Using the exact sequence $0 \rt R(I,M) \rt M[t] \rt E \rt 0$ and the
corresponding long exact sequence for local cohomology  we get  (a) and (b).

Since $ H^{i}_{\M}\left(R(I,M)\right)$ is *-Artinian for all $i \geq 0$ we get
$ H^{i}_{\M}\left( E\right)$ is *-Artinian for   $0 \leq i \leq l -1$. So
$ H^{i}_{\M}\left( L^I(M)\right)$ is *-Artinian for   $0 \leq i \leq l -1$.
\end{proof}

\s \label{perexam2}
Let $I = (x_1,\ldots,x_{m})$.
 Set $S = A[X_1,\ldots,X_m]$.  We have a
surjective  homogeneous
homomorphism of $A$-algebras, namely $ \phi: S \rt R(I)$ where $\phi(X_i)
= x_it$. We also have the natural map $\psi: R(I) \rt G_I(A)$.
Set $G = G_I(A)$ and $R = R(I)$.
Note that
\[
 \phi(\M_S) = \M_{R}, \quad
 \psi(\M_R) = \M_{G}\quad \text{and} \quad \psi\circ\phi(\M_S) = \M_G.
\]
By graded independence theorem cf. \cite[13.1.6]{BSh} it does not matter which ring we
 use to compute local cohomolgy.

\s \label{basechangeC}\textbf{Local cohomology and base change}
Let $\phi \colon (A,\m) \rt (A',\m')$ be the local ring homomorphism as discussed in
\ref{AtoA'} i.e., either (i) $A'$ is a quotient of $A$ or (ii)  a flat $A$-algebra with $\m' = \m A'$.
We use notation as in \ref{perexam2}. Set $S' = A'[X_1,\ldots,X_m]$.  Notice
$S' = S \otimes_A A'$ and
\begin{alignat*}{2}
L^I(M)\otimes_S S' &=  L^I(M) \otimes_A A' &  &= L^{I'}(M') \\
 G_I(M)\otimes_S S' & =  G_I(M)\otimes_A A' &  &=  G_{I'}(M')
\end{alignat*}
In case (i) by graded independence theorem cf. \cite[13.1.6]{BSh} it does not matter whether
we compute local cohomolgy \wrt \ $S$ or $S'$.
In case (ii) by graded flat base theorem it follows that for all $i \geq 0$ we have
\begin{alignat*}{2}
H^i(L^{I'}(M)) &= H^i(L^I(M))\otimes_S  S'& &= H^i(L^I(M))\otimes_A A' \\
H^i(G_{I'}(M)) &= H^i(G_I(M))\otimes_S  S'& &=  H^i(G_I(M))\otimes_A A'.
\end{alignat*}

We  compute the zeroth local cohomolgy of $L^{I}(M)$ \wrt  \ $\M = \M_{R(I)}$ and $R(I)_+$.
It is convenient to define the following $R(I)$-module
\[
\R^{I}(M) = \bigoplus_{i \geq 0}\frac{\wt{I^{i+1}M}}{ I^{i+1}M} \quad = \bigoplus_{i = 0}^{\rho^{I}(M) - 1}\frac{\wt{I^{i+1}M}}{ I^{i+1}M}
\]

\begin{proposition}
\label{0lc}
(Hypothesis as in \ref{setupgen}.) We have
\begin{enumerate}[\rm 1.]
\item
$\displaystyle{H^{0}_{R_+}(L^I(M)) = \R^{I}(M).}$
\item
If $I$ is an ideal of definition for $M$ then $H^{0}_{\M}(L^I(M)) = \R^{I}(M)$.
\end{enumerate}
\end{proposition}
\begin{proof}
Let $I = (x_1,\ldots,x_m)$. Set $S = A[X_1,\ldots,X_m]$. As
described in \ref{perexam2} consider $R$ as an $S$-module.
Set $\aF = S_+$ and $\n = \M_S$.
Since $\grade(I,M) > 0$ we have $\wt{I^nM} = I^nM$ for all $n \gg 0$.
So $\R^{I}(M)_n = 0$ for all $n \gg 0$. Therefore $\R^{I}(M) \subseteq \Gamma_{\aF}(L^I(M))$.
We claim that
$\Gamma_{\aF}(L^I(M)) \subseteq \R^{I}(M).$
%\end{equation*}
Let $\ov{p} \in \Gamma_{\aF}(L^I(M))$ be homogeneous of degree $i$. Note that
 $\aF^{l}(\ov{p})= 0$ implies that $I^lp \subseteq I^{l+i+1}M$.
So $p \in ( I^{l+i+1}M \colon_M I^l) \subseteq \wt{I^{i+1}M}$.
Therefore $\ov{p} \in (\R^{I}(M))_i$. So
$\Gamma_{\aF}(L^I(M)) \subseteq \R^{I}(M)$.
Thus $\Gamma_{\aF}(L^I(M)) = \R^{I}(M)$.

Note that in general
 $ \Gamma_{\n}(L^I(M)) \subseteq \Gamma_{\aF}(L^I(M))$. So
 $ \Gamma_{\n}(L^I(M)) \subseteq \R^{I}(M)$.
If $I$ is an ideal of definition for $M$ then by \ref{basechangeC} we may assume $I$ is $\m$-primary. Set
$\q = (I, X_1,\ldots,X_s)R = IR + \A$. Then since $\sqrt{\q} = \n$ we have that
$\Gamma_{\q}(L^I(M)) = \Gamma_{\n}(L^I(M))$. We prove $\R^{I}(M) \subseteq \Gamma_{\q}(L^I(M))$.

\textbf{Claim:} $I^s\A^r \R^I(M) = {0}$ if $s+r \geq \rho^I(M)$.

\noindent If $\ov{p} \in \wt{I^{j+1}M}/I^{j+1}M$ then $\A^r \ov{p} \subseteq \wt{I^{j+r+1}M}/I^{j+r+1}M$.
Notice
\[
I^s\A^r\ov{p} \subseteq \frac{\wt{I^{j+r+s+1}M} + I^{j+1}M}{I^{j+1}M}.
\]
Since $r+s \geq \rho^I(M)$ we have $\wt{I^{j+r+s+1}M} = I^{j+r+s+1}M \subseteq I^{j+1}M$.
So $I^s\A^r\ov{p} = {0}$.

Notice $\q^{l}= \sum_{r+s = l}I^s\A^rR$.
Consequently  $\q^{l}\R^{I}(M) = 0$ if $l \geq \rho^{I}(M)$.
\end{proof}

\section{The first fundamental exact sequence and  applications}
\s
In this section the setup is as in \ref{setupgen}.

\s
The first fundamental exact sequence is
\begin{equation}
\label{dag}
0 \xar G_{I}(M) \xar L^I(M) \xrightarrow{\Pi} L^I(M)(-1) \xar 0.
\end{equation}

Using (\ref{dag}) and Proposition \ref{0lc} we immediately get
%yyyy
\begin{theorem}
\label{mthm}
 (with hypothesis as in \ref{setupgen}.) We have
\[
 \rho^{I}(M) \leq \max\{0, a_{1}\left(G_{I}(M)\right) + 1 \}.
\]
\end{theorem}
\begin{proof}
Set $R = R(I)$, $\A = R(I)_+$ and $a_1 = a_1(G_I(M))$.
We  take local
cohomology \wrt \ $\A$. Set $b_0 = \e(H^{0}_{\A}(L^{I}(M))$. Since
$\R^{I}(M) = H^{0}_{\A}(L^{I}(M))$ we have
$ \rho^{I}(M) = \max \{0, b_0 + 2 \}$.
 Using (\ref{dag}) we have an exact sequence
\begin{align*}
0  &\xar H^{0}_{\A}(G_{I}(M)) \xar H^{0}_{\A}(L^{I}(M)) \xar H^{0}_{\A}(L^{I}(M))(-1) \\
 &\xar H^{1}_{\A}(G_{I}(M))
\end{align*}
Therefore $b_0 \leq a_1 -1$. This establishes the assertion of the Theorem.
\end{proof}
%yyy
The following Lemma is well-known.
\begin{lemma}
\label{monart} Let $R =\bigoplus_{i \geq 0}R_i $ be a graded ring.
 Let $L$ be a graded *-Artinian
$R$-module.Then
\begin{enumerate}[\rm 1.]
\item
$L_n = 0$ for all $n \gg 0$.
\item
If $\psi \colon L(-1) \rt L$ is a monomorphism then $L = 0$.
\item
If $\phi \colon L \rt L(-1)$ is a monomorphism then $L = 0$.
\end{enumerate}
\qed
\end{lemma}

%yyyy
\begin{proposition}
\label{eqndepth}
 (with hypothesis as in \ref{setupgen}.)
Set $R = R(I)$ and  $\M = \M_R$ the $*$-maximal ideal of $R$.
Then
\begin{enumerate}[\rm 1.]
\item
If $i \leq \depth M - 1$ and $H^{i}_{\M}(G_I(A)) = 0$ then $ H^{i}(L^{I}(M)) = 0.$
\item
Let $s \leq \depth M -1$. We have
$$ H^{i}(L^{I}(M)) = 0  \ \text{for $i = 0,\ldots,s$ \ \text{iff}  \ } \
  H^{i}(G_I(M)) = 0 \ \text{for $i = 0,\ldots,s$.} $$
\end{enumerate}
\end{proposition}
\begin{proof}
If $\depth M = 0$ there is nothing to prove. So assume $\depth M > 0$.
Set $G = G_I(M)$ and $L = L^I(M)$.
We use (\ref{dag}) and the corresponding long exact sequence in local cohomolgy.

1. If $H^{i}_{\M}(G) = 0 $  then we have an injective map
$ H^{i}_{\M}(L) \xar
H^{i}_{\M}(L)(-1)$. Since $i < \depth M - 1$,  $ H^{i}_{\M}(L)$ is
*-Artinian.
Using Lemma \ref{monart}
we get
 $ H^{i}_{\M}(L) = 0$.

2. Using
 (\ref{dag}) and the corresponding long exact sequence in local cohomolgy,
it
follows that if $ H^{i}_{\M}(L^{I}(M)) = 0 $ for $0 \leq i \leq s$
 then $H^{i}_{\M}(G_I(M)) = 0 $ for $0 \leq i \leq s$.
 The converse follows from 1.
\end{proof}

An easy application of previous Proposition is the following extension to Rees modules a  result of Huckaba and Marley regarding depth of
the Rees ring \cite{HMr}.
\begin{proposition}
\label{HMR}
 (with hypothesis as in \ref{setupgen}.)
If $\depth G_I(M) < \depth M$ then $\depth R(I,M) =  \depth G_I(M) + 1$.
\end{proposition}
\begin{proof}
Set $R = R(I)$, $\M = \M_{R(I)}$ and
  $s = \depth G_I(M)$.
Using Proposition \ref{eqndepth} we get that $ H^{i}_{\M}(L^{I}(M))$ is zero
 for $i<s$.
Using
 (\ref{dag}) and the corresponding long exact sequence in local cohomolgy,
it
follows that $ H^{s}_{\M}(L^{I}(M)) \neq 0$.
Note that $ H^{0}_{\M}(R(I,M)) = 0$.  Since $s < \depth M$, it follows
from Proposition \ref{fgenExt}  that
 $ H^{i}_{\M}(R(I,M))$ is zero for $1 \leq i<s+1$
and is non-zero for $i =s +1$.
So $\depth R(I,M) =  s + 1$.
\end{proof}
\s \label{binvarNew} \textbf{A new invariant:} For $i = 0,\ldots, \depth M -1$,
the modules $H^{i}(L^I(M))$  are *-Artinian. This enables us to define a
new invariants
\begin{align*}
b^{I}_{i}(M) &= \e\left(H^{I}(L^I(M))\right) \quad i = 0,\ldots,\depth M -1; \\
b^I(M) &= \max \{ b^{I}_{i}(M) + i \mid 0 \leq i \leq \depth M -1 \}.
\end{align*}
Since $H^i(G_I(M))$ is *-Artinian we can define
\begin{align}
a_{i}^{*}(G_I(M)) &= \e \left(H^{i}(G_I(M))\right) \quad \text{for $i \geq 0$}; \\
\label{reg*}\reg^*(G_I(M)) &= \max\{ a_{i}^{*}(G_I(M)) + i \mid i = 0,\ldots,\dim M \}.
\end{align}
Using (\ref{dag}) we get
\begin{equation}
\label{barelation}
b^{I}_{i}(M) \leq a_{i+1}^{*}(G_I(M)) -1.
\end{equation}
Consequently we get
\begin{equation}
\label{breg}
b^I(M) \leq \max \{0, \reg(G_I(M)) -2\}.
\end{equation}
Recall $$a_i(G_I(M)) = \e \left(H^{i}_{G_I(A)_+}(G_I(M))\right).$$
\begin{remark}
If $I$ is an ideal of definition for $M$ then   $a_{i}^{*}(-) = a_i(-)$ and
so $\reg^*(-) = \reg(-)$.
\end{remark}

\s\label{marextnD}
Let $I$ be an  ideal in a  local ring $A$. Set $G = G_I(A)$ and let $s = \depth G_I(A)$.
If $I$ is $\m$-primary and $A$ is \CM \ Marley \cite[Theorem 2.1]{Mar}  shows
$a_{s}(G) < a_{s+1}(G)$.
  This was generalized to an arbitrary ideal $I$ with $s \leq \grade I - 1$ by Hoa
\cite[Theorem 5.2]{Hoa2}. See also  \cite[Proposition 6.1]{Tr} for a different proof.
The same proof goes through for modules.  An easy consequence of our investigations is the following:
\begin{corollary}\label{marextn}
 (with hypothesis as in \ref{setupgen}.)
 Set $s = \depth G_I(M)$ and assume $s \leq \grade(I,M) - 1$. Then $a_{s}^{*}(G_I(M)) <
a_{s+1}^{*}(G_I(M))$.
\end{corollary}
\begin{remark}
If $M =A$ and
if $I$ is \emph{not} $\m$-primary then this result is different from Hoa's generalization.
\end{remark}
\begin{proof}[Proof of Corollary \ref{marextn}]
For $i \geq 0$,
set $a_{i}^{*} = a_{i}^{*}(G_I(M))$ and $b_i = b_{i}^{I}(M)$.
By Proposition \ref{eqndepth} we get that
\[
H^i(G_I(M)) = H^i(L^I(M)) = 0 \quad \text{for} \ i < s.
\]
Using
 (\ref{dag}) and the corresponding long exact sequence in local cohomolgy
 we get $a_{s}^{*} \leq b_{s}$. Also by
 (\ref{barelation}) we get $b_{s} \leq a_{s+1}^{*} - 1$. So we have
 $a_{s}^{*} < a_{s+1}^{*}$.
\end{proof}

\section{Ratliff-Rush Filtration and superficial elements}
Ratliff-Rush filtration does not behave well \wrt \ superficial
elements; see \cite[Section 2]{SwR}. However  we construct two 
exact sequences which arise naturally in this context.

\s\label{setuprr} Throughout this section the setup is as in \ref{setupgen}.
\s \label{4es}
 Let $x$ be $M$-superficial \wrt \  $I$.
We consider the homomorphism:
\[
 \alpha_{n-1}^{x} : \frac{\wt{I^nM}}{I^nM} \xar  \frac{\wt{I^{n+1}M}}{I^{n+1}M} \quad \text{defined by}\ \  \alpha_{n-1}^{x}(u +I^nM) = xu + I^{n+1}M.
\]

This yields the    exact sequence:
\begin{equation}
\label{supex}
 0 \xar \frac{(I^{n+1}M\colon_M x)}{I^nM} \xar \frac{\wt{I^nM}}{I^nM} \xrightarrow{\alpha_{n-1}^{x}}\frac{\wt{I^{n+1}M}}{I^{n+1}M}.
\end{equation}

If $x$ is $M$-superficial \wrt \ $I$ we have $(I^{n+1}M \colon_M x) = I^nM$ for all $n \gg 0$. Set
$$\rho^{I}(x,M) = \min\{i \mid (I^{n+1}M \colon_M x) = I^nM \ \text{for all} \ n \geq i \}.$$
An  application of (\ref{supex}) is  the following
\begin{corollary}
(with the hypothesis as above) $\rho^{I}(x,M) = \rho^{I}(M)$. Thus
$\rho^{I}(x,M)$ is independent of superficial elements.
\end{corollary}
\begin{proof}\begin{align*}
\text{Notice} \quad (I^{n+1}M\colon_M x) &= I^nM  \ \text{for } \ n \geq \rho^{I}(M) \ \text{and}\\
 &= \wt{I^nM} \left(\neq I^nM\right) \  \text{for } \ n =  \rho^{I}(M) - 1.
\end{align*}
\end{proof}

\begin{remark}
\label{posdepth}
If $\depth G_I(M) > 0$ and $x$ is $M$-superficial \wrt \ $I$ then $I^{n+1}M \colon_M x = I^nM$ for
each $n \geq 1$. So $\rho^{I}(M) = 0$. Conversely if $\rho^I(M) = 0$ i.e., $\wt{I^nM} = I^nM$ for all $n \geq 1$, then we get
$(I^{n+1}M\colon_M x) = I^nM  \ \text{for all } \ n \geq 0$. So $x^*$ is $G_I(M)$-regular.
\end{remark}

\begin{remark}
Let $M$ be an $A$-module with
 $\grade(I,M) \geq 2$. Let $x$ be an
$M$-superficial element \wrt \ $I$ and set $N = M/xM$.  Clearly
$\ov{\wt{I^nM}} \subseteq \ \wt{I^nN} $ for all $n \geq 0$.  We
 have the following
exact sequence
\begin{equation}
\label{supexN}
 0 \xar \frac{(I^{n+1}M\colon_M x)}{I^nM} \xar \frac{\wt{I^nM}}{I^nM} \xrightarrow{\alpha_{n-1}^{x}}\frac{\wt{I^{n+1}M}}{I^{n+1}M} \xrightarrow{\rho_n} \frac{\wt{I^{n+1}N}}{I^{n+1}N};
\end{equation}
Here $\rho_n$ is the natural quotient map (defined since $\ov{\wt{I^nM}} \subseteq \ \wt{I^nN} $ for all $n \geq 0$.)
\begin{equation}
\label{supexN0}
\text{In particular} \quad \rho_0 \colon\frac{\wt{IM}}{IM} \xar \frac{\wt{IN}}{IN}  \quad \text{is injective.}
\end{equation}
\end{remark}
Another motivation for this paper  was to determine whether
 (\ref{supexN}) is part of a longer exact sequence. A surprising  application of
(\ref{supexN}) and (\ref{supexN0}) is the following proof of Sally-Descent
For $M = A$
see \cite{Sa2b}, \cite[2.2]{HM}. The general case was proved by the author \cite[8.(2)]{Pu1}.

\begin{theorem}
\label{sd}
Let $(A,\m)$ be a local ring, $M$ be a finite $A$-module with $\depth M \geq r + 1$ and
 $I$ be an ideal of definition for $M$. Let  $x_1,\ldots,x_r$ be an $M$superficial sequence
\wrt \ $I$.
Set $N = M/(x_1,\ldots,x_r)M$. If $\depth G(N) \geq 1 $ then $\depth G(M) \geq r+1$.
\end{theorem}
\begin{proof} In view of \cite[Theorem 8(1)]{Pu1}
 it suffices to consider the
case when $r = 1$. We may also assume residue field of $A$ is infinite ( see \ref{AtoA'}). Let $x$ be $M$-superficial \wrt \ $I$.  If $\depth G(N) > 0$ then by Remark \ref{posdepth}
we get $\wt{I^nN}= I^n N$ for all $n \geq 1$. So by (\ref{supexN}) we get
$\wt{IM} = IM$. Using (\ref{supexN0}) recursively we get
$\wt{I^nM}= I^n M$ for all $n \geq 1$ and
 $(I^{n+1}M \colon_M x) = I^nM $ for all
$n\geq 1$. By Remark \ref{posdepth} we have $x^*$ is $G_I(M)$ regular and
by \cite[7]{Pu1} we have $G_I(M)/x^*G_I(M) \cong G_I(N)$. It follows
 that $\depth G_I(M) \geq 2$.
\end{proof}

\section{The second fundamental exact sequence and applications}

\s \textbf{Setup:}\label{setup}

 $(A,\m)$ is local and  $M$ is an $A$-module.
The ideal $I = (x_1,\ldots,x_m)$ in $A$ is
 an \emph{ideal of definition of $M$}.
We assume  $\depth M \geq 2$.
 Set  $S = A[X_1,\ldots,X_m]$ and consider
$R(I)$ as a quotient ring of $S$ (see \ref{perexam2}). Set  $\M = \M_S$. Throughout
we take local cohomolgy \wrt \ $\M$. We denote
$H^{i}_{\M}(-) $ by $H^i(-)$.

\s \textbf{The second fundamental exact sequence:}
Let $x $ be $M$-superficial \wrt \  $I$. There exists  $c_i \in A$ such that $x
= \sum_{i= 1}^{s}c_ix_i$. Set $X = \sum_{i= s}^{l}c_iX_i$ and  $N
= M/xM$. We  have an exact sequence:
\begin{equation}
\label{dagg}
0 \xar \B^{I}(x,M) \xar L^{I}(M)(-1)\xrightarrow{\Psi_X} L^{I}(M) \xrightarrow{\rho}  L^{I}(N)\xar 0,
\end{equation}
where $\Psi_X$ is left multiplication by $X$ and
\[
\B^{I}(x,M) = \bigoplus_{n \geq 0}\frac{(I^{n+1}M\colon_M x)}{I^nM}.
\]

\s \label{longHpara} The  exact sequence below connects the local cohomolgy of $L^I(M)$ and
$L^I(N)$.
\begin{equation}
\label{longH}
\begin{split}
0 \xar \B^{I}(x,M) &\xar H^{0}(L^{I}(M))(-1) \xar H^{0}(L^{I}(M)) \xar H^{0}(L^{I}(N)) \\
                  &\xar H^{1}(L^{I}(M))(-1) \xar H^{1}(L^{I}(M)) \xar H^{1}(L^{I}(N)) \\
                 & \cdots \cdots \\
                &\xar H^{i}(L^{I}(M))(-1) \xar H^{i}(L^{I}(M)) \xar H^{i}(L^{I}(N))\cdots \\
\end{split}
\end{equation}
To see this break  (\ref{dagg}) into two short exact
sequences.
  The result follows since $\B^{I}(t,M)$ has finite length.
Notice if we take $n$-th degree of top row of (\ref{longH}) we recover
(\ref{supexN}).

Next we prove a crucial theorem of this paper.

\begin{theorem}\label{artinmb}
Let $(A,\m)$ be a Noetherian local ring, $M$  a finitely generated $A$-module
and  $I$  an  ideal. For $0 \leq i \leq  \depth M - 1$ we have
\begin{enumerate}[\rm 1.]
\item
The modules $H^{i}(L^I(M))$ are  *-Artinian.
\item
If $I$ is an ideal of definition for $M$ then
\begin{enumerate}[\rm a.]
\item
 $H^{i}(L^I(M))_n$  has finite length
for all $n \in \mathbb{Z}$.
\item
$\lambda(H^{i}(L^I(M))_n)$  coincides with a polynomial for all $n \ll 0$.
\end{enumerate}
\end{enumerate}
\end{theorem}
\begin{proof}
 Proposition \ref{fgenExt} implies 1. We prove 2.a. by induction on $\depth M$.
By  \ref{basechangeC} we may assume that the residue field of $A$ is infinite.
When $\depth M  = 1$, Proposition \ref{0lc} implies the result.
Assume the result for modules with depth $l$. We prove when
$\depth M = l+1$. Let $x$ be $M$-superficial \wrt \ $I$.
Set $N = M/xM$. By 1. the $R(I)$-modules $H^{i}(L^{I}(M))$ are *-Artinian
for $0 \leq i \leq  l$. So    $b_i = \e(H^i(L^I(M))) <
\infty $. By induction hypothesis,  for $0 \leq i \leq  l -1$, the $A$-module $H^i(L^I(N))_j$
have finite length for all $j \in \Z$. Fix
$i$ with $1\leq i \leq l$. We prove by induction on $m$ that $H^i(L^I(M))_{b_i
- m}$ has finite length for all $m \geq 0$.
 When $m = 0$ we use a part
of  (\ref{longH})
\[
H^{i-1}(L^{I}(N))_{b_i + 1} \xar H^{i}(L^{I}(M))_{b_i }\xar H^{i}(L^{I}(M))_{b_i + 1} = 0.
\]
So $ H^{i}(L^{I}(M))_{b_i }$ has finite length. Assume the result for $m = c$. We prove it for $m = c+1$.
We use a part of  (\ref{longH})
\[
H^{i-1}(L^{I}(N))_{b_i - c} \xar H^{i}(L^{I}(M))_{b_i -c -1 }\xar H^{i}(L^{I}(M))_{b_i -c}.
\]
Since $H^{i-1}(L^{I}(N))_{b_i - c}$ and $H^{i}(L^{I}(M))_{b_i -c}$ have finite length it follows that

\noindent  $H^{i}(L^{I}(M))_{b_i -c-1}$ has finite length.
Thus we have shown that $H^i(L^I(M))_j$ has finite length for all $j \in \Z$ and
$i = 1,\ldots,l-1$. The case $i = 0$ is taken care by Proposition \ref{0lc}.

2.b. By \ref{basechangeC} we may assume $I$ is $\m$-primary. Set $G = G_I(A)$ and $L = L^I(M)$.
Note that $H^{i}(G_I(M))$ is *-Artinian $G$-module for all $i \geq 0$.
We use (\ref{dag}) and the corresponding long exact sequence in cohomolgy to get
\begin{equation*}
 H^{i}(G_I(M)) \xrightarrow{u^i} H^{i}(L) \rt H^{i}(L(-1)) \xrightarrow{\delta^i} H^{i+1}(G_I(M)). \tag{$\dagger$}
\end{equation*}
% \tag{\dagger}
%\end{equation*}
Set
$D^i = \image(u^i) \quad \text{and} \quad E^i = \image(\delta^i).$
Clearly $D^i, E^i$ are *-Artinian $G$-modules and therefore  by \ref{polyartin},  $\lambda(D^{i}_{n})$ and
$\lambda(E^{i}_{n})$ are polynomials for all $n \ll 0$.
  Using ($\dagger$) we get
\[
\lambda\left(H^{i}(L)_n \right) -  \lambda\left(H^{i}(L)_{n-1}\right) = \lambda\left(D^{i}_{n}\right) -  \lambda\left(E^{i}_{n}\right)
\]
for $i = 0,\ldots, \depth M - 1$.
It follows that $\lambda\left(H^{i}(L)_n \right)$ is  polynomial
for $n \ll 0$.
\end{proof}

\s \label{binvar}
 If $x$ is $M$-superficial \wrt \ $I$ then by
(\ref{longH}) it follows   that
\begin{equation}
\label{bmodx}
b_{i}^{I}(M/xM) \leq \max\{b_{i}^{I}(M), b_{i+1}^{I}(M) +1 \} \quad  \text{for $0 \leq i \leq \depth M -2$}.
\end{equation}
So we have $b^{I}(M/xM) \leq b^{I}(M)$.

%\label{verEs}

\section{Powers of $I$}
\s In this section the setup is as given in \ref{setup}.

It is easy to see that $R(I^l) =R(I)^{<l>} $.
We also have
\begin{equation}
\label{ver}
L^{I^{l}}(M)(-1)\  = \  \bigoplus_{n \geq 0}\frac{M}{I^{nl}M} \  = \ L^{I}(M)(-1)^{<l>}.
\end{equation}
\textbf{Caution:} It can be directly seen that $G_I(A)^{l} \neq G_{I^l}(M)$.

\s \label{h1asymm}
 Recall that local cohomology commutes with the Veronese functor cf. \ref{verEs}.
In particular  if  $H^{1}(L^I(M)) = 0$ then $\depth G_{I^n}(M) \geq 2$ for all $n \gg 0$;
(Use (\ref{ver}), Proposition \ref{0lc} and Proposition \ref{eqndepth}).

As an  application of the previous deliberations  we give another proof of a result due to
Huckaba and Huneke \cite[3.8]{HHu}.
\begin{theorem}
\label{mynor}
Let $(A,\m)$ be a Noetherian local ring with depth $s \geq 2$. Let $I$ be an $\m$-primary
asymptotically normal ideal. Then $\depth G_{I^n}(A) \geq 2$ for all $n \gg 0$.
If $I$ is normal then $\depth G_{I^n}(A) \geq 2$ for all $n \geq  \max \{1, \reg(G_I(A)\}$.
\end{theorem}
\begin{proof}
Set $ u = \max \{ b_{0}^{I}(A), b^{I}_{1}(A)\} + 2$.
It can be easily verified that for $i = 0, 1$ and $l \geq u$
\begin{equation*}
\e\left(H^{i}(L^{I^l}(A)(-1)) \right) = \e\left(H^{i}(L^{I}(A)(-1)^{<l>}) \right) \leq 0. \tag{*}
\end{equation*}
Set $v = \min \{ m \mid I^n \ \text{is integrally closed for all} \ n \geq m \}$.
Fix $l \geq \max \{ u,v \}.$ Set $K = I^l$. Using \ref{AtoA'}(ii).b  we may assume there exists $x \in K$
such that the $B = A/(x)$ ideal $J = K/(x)$ is integrally closed.  Using
(*) we get that
\begin{align*}
\e(H^i(L^K(A)) &\leq -1 \quad \text{for $i = 0,1$ and so by (\ref{longH})} \\
\e(H^0(L^K(B)) &\leq 0.
\end{align*}
Notice $ H^0(L^K(B))_0 = \wt{J}/J = 0$ since $J$ is integrally closed. Thus $H^0(L^K(B)) = 0$.
By Proposition \ref{eqndepth}.2  we get that $\depth G_K(B) \geq 1$ and so by Sally descent we get
$\depth G_K(A) \geq 2$.

If $I$ is normal then  $v = 1$. Also by (\ref{breg}) we have $u \leq \max\{ \reg(G_I(A)),0 \}$.
Thus if  $l \geq \max \{1, \reg(G_I(A))\}$ then $\depth G_{I^l}(A) \geq 2$.
\end{proof}

\begin{definition}
Let $M$ be an \ $A$-module of depth $s \geq 1$. Set $L = L^I(M)$ for
$i = 0, \ldots, s-1$.
\begin{align*}
\xi_I(M) &:= \underset{0 \leq i \leq s-1}\min\{ \ i \ \mid H^{i}(L)_{-1} \neq 0 \ \text{or} \ \lambda(H^{i}(L)) = \infty \}.\\
\amp_I(M) &:= \max\{\ |n| \ \mid H^{i}(L)_{n-1} \neq 0 \  \text{for} \ i = 0,\ldots, \xi_I(M) - 1 \}.
\end{align*}
\end{definition}
Since $\lambda(H^0(L)$ is finite and since $  H^{0}(L)_{-1}  =0 $ we
get $\xi_I(M) \geq 1$.
The reason for calling the first constant above, $\xi_I(M)$, is the following theorem.
\begin{theorem}
\label{asyym}
Let $(A,\m)$ be local and $M$ a   $A$-module of depth $s \geq 1$. We have
\[
 \depth G_{I^l}(M) = \xi_I(M) \ \text{for all} \ l > \amp_I(M).
\]
\end{theorem}
As a corollary we immediately get
\begin{corollary}
\label{asyymCor}
Let $(A,\m)$ be local and $M$ an   $A$-module. We have
$ \depth G_{I^n}(M) =$ constant for all $n \gg 0$.
\end{corollary}
\begin{proof}
Notice
\[
l = \depth M \geq \underset{n \geq 1}\max\{ \depth G_{I^n}(M) \}.
\]
If $l = 0$ then $\depth  G_{I^n}(M) = 0$ for all $n \geq 1$.  If $l  \geq 1$ the result
follows from Theorem \ref{asyym}.
\end{proof}

\begin{proof}[Proof of Theorem \ref{asyym}]
Set $E^i = H^{i}\big(L^I(M)(-1)\big)$ for
$i = 0, \ldots, r-1$. Set $u = \xi_I(M)$. Notice  Fix $l > \amp_I(M)$.

For $i = 0, \ldots, u-1$ notice  $(E^{i})^{<l>} = 0$.
So we have $H^{i}\big(L^{I^l}(M)(-1)\big) = 0$ for $i = 0, \ldots, u-1$.
Therefore by Proposition \ref{eqndepth} we have
$\depth G_{I^l}(M) \geq  u$.  

Suppose if possible $\depth G_{I^l}(M) > u $. Note $\depth M > u$. Since

\noindent 
$H^{u}\big(L^{I^l}(M)(-1)\big) = 0$ we have $E^{u}_{0} = 0$ and $E^{u}_{nl} = 0$
for all $n  \in \Z$. By Theorem \ref{artinmb}, $\lambda(E^{u}_{n})$ is polynomial
for all $n \ll 0$. Since $E^{s}_{nl} = 0$ for all $n \leq - 1$, we get
$E^{u}_{n} = 0$ for all $n \ll 0$. As $E^u$ is *-Artinian, $E^{u}_{n} = 0$ for all $n \gg 0$.
Therefore $\lambda(E^u)$ is finite.
This implies $\xi_I(M) \geq u+1$, a contradiction.
Thus $\depth G_{I^l}(M) = u$.
\end{proof}

\begin{remark}
\label{asympconse}
If for some $l \geq 1$ we have $\depth G_{I^l}(M) = u$ then $\depth G_{I^{nl}}(M) \geq u$
for all $n \geq 1$. Thus $\xi_I(M) \geq u$. In particular
if $ G_{I^l}(M) $
is \CM \ for some $l$ then $\xi_I(M) = r$ i.e. $G_{I^n}(M)$ is \CM \ for all $n \gg 0$.
\end{remark}

\begin{proposition}
\label{pgcm}
Let $(A, \m)$ be local,
Let $M$ be a Cohen-Macaulay $A$-module of dimension $r$ and let $I$ be an
ideal of definition for $M$.
 If \ $G_{I^n}(M)$ is \CM \  for some $n \geq 1$ then $G_{I}(M)$ is generalized \CM.
\end{proposition}
\begin{proof}
By Remark \ref{asympconse}
we get $\xi_I(M) = r$.
Therefore  $H^{i}( L^{I}(M))$ has finite length  for
all  $i = 0,1,\ldots,r-1$. Using (\ref{dag}) and the corresponding long exact  sequence
in cohomolgy,
 we get that $H^{i}( G_I(M))$ has finite length for
all  $i = 0,1,\ldots,r-1$. So $G_I(M)$ is generalized \CM.
\end{proof}

We first prove the following general result.
\begin{proposition}
\label{propxi1}
(with hypothesis as in \ref{setup})
 Let $x \in I$ be $M$-superficial \wrt \ $I$. Set
$N = M/xM$. If $\wt{IM} = IM$ and $\wt{IN} \neq IN$ then $\xi_I(M) = 1$.
\end{proposition}
\begin{proof}
Using (\ref{longH}) it follows that $H^1(L^I(M)(-1))_0 \neq 0$. So by
definition we get $\xi_I(M) \leq 1$. However $\xi_I(M) \geq 1$ always. Therefore
 $\xi_I(M) = 1$.
\end{proof}

This proposition above gave me an idea for the following example of a  \CM  \ local ring $A$ of dimension $d =  r +s$ with $r \geq 1, s\geq 1$ and $I$ an $\m$-primary ideal 
 with $\xi_I(A) = r$.
\begin{example}
\label{arbitxi}
Let $(A_0,\m_0)$ be a \CM \ local ring of dimension $s \geq 1$. Let $l \geq 0 $  and let $\A_0$
 be an $\m_0$-primary ideal with
$\wt{\A_{0}^{m}} \neq \A_{0}^{m} $ for $m =  1,\ldots, l +1$. Let $r$ be an integer with $1 \leq r \leq l$.
Set $A_r = A_0[X_1,\ldots,X_r]_{(\m_0,X_1,\ldots,X_r)}$ and $\m_r = (\m_0,X_1,\ldots,X_r)$.
Note that $A_r$ is \CM \ of dimension $d = r +s$ and $\A_r = (\A_0,X_1,\ldots,X_r)$ is
$\m_r$-primary. Notice $G_{\A_r}(A_r) = G_{\A_0}(A_0)[X_{1}^{*},\ldots,X_{r}^{*}]$.
 We claim
$\xi_{\A_r}(A_r) = r$.

\noindent \textit{Proof of the claim:}
Notice that for $i = 0,\ldots,r-1$, the ring $A_i$ is a quotient of $A_r$ and
 $\A_rA_i = \A_i$.
Set $L(A_i) = L^{\A_i}(A_i)  = L^{\A_r}(A_i)$. If $i \geq 1$ then $X_{1}^{*},\ldots,X_{i}^{*}$
is a $G_{\A_r}(A_i)$ regular sequence and so by Proposition \ref{eqndepth} we get
\begin{equation}
\label{arbit1}
  H^{j}(L(A_i)) = 0 \quad \text{for} \ 0\leq j \leq i-1.
\end{equation}
 We prove for $i \geq 0$,
\begin{equation}
\label{arbit2}
  H^{i}(L(A_i))_{n} \neq 0 \quad \text{when} \ -i\leq n \leq l-i.
\end{equation}
If we prove (\ref{arbit2}) then we have $H^{r}(L(A_r))_{-1} \neq 0$. So by (\ref{arbit1})
and Theorem \ref{asyym} we get $\xi_{\A_r}(A_r) = r$.
We prove (\ref{arbit2}) by induction on $i$.
If $i = 0$ then $H^0(L(A_0))_{n} = \wt{\A_{0}^{n+1}}/\A_{0}^{n+1} \neq 0$ for $0 \leq n \leq l$
by hypothesis on $\A_0$.  We assume (\ref{arbit2}) if $i = p-1 \geq 0$ and prove it for
$i = p$.
Using (\ref{arbit1}) and (\ref{longH}) we get an exact sequence
\begin{equation}\label{arbit3}
0 \xar H^{p-1}(L(A_{p-1})) \xar  H^{p}(L(A_p))(-1).
\end{equation}
By induction hypothesis $H^{p-1}(L(A_{p-1}))_n \neq 0$ for $-(p-1) \leq n \leq l -(p-1)$. Using
(\ref{arbit3}) we get $H^{p}(L(A_p))_n \neq 0$ for $ -p \leq n \leq l -p$. So (\ref{arbit2}) holds for $i =p$.
This proves (\ref{arbit2}) and as discussed before it implies $\xi_{\A_r}(A_r) = r$.
\end{example}

An easy corollary of Proposition \ref{propxi1} is the following
\begin{proposition}\label{curiousm}
Let $(A,\m)$ be a Noetherian local ring with $\depth A \geq 2$. Let $I$ be an $\m$-primary
 Ratliff-Rush closed ideal and assume that $\depth G_{I^n}(A) \geq 2$. If
$x$ is $A$-superficial \wrt \ $I$ then $I/(x)$ is a Ratliff-Rush closed ideal.\qed
\end{proposition}

Using Huneke and Huckaba's result on normal ideals cf. Theorem\ref{mynor} and Proposition \ref{curiousm} we immediately get:
\begin{corollary}
\label{curiousnormal}
Let $(A,\m)$ be a Noetherian local ring with $\depth A \geq 2$. Let $I$ be an $\m$-primary
 Ratliff-Rush closed ideal and assume that $I$ is asymptotically normal. If
$x$ is $A$-superficial \wrt \ $I$ then $I/(x)$ is a Ratliff-Rush closed ideal.\qed
\end{corollary}

Next we give an example of a \CM \ local ring $(A,\m)$ of dimension 2 such that $G_{\m}(A)$ is not generalised \CM. In view of Proposition \ref{pgcm} this implies $\xi_{\m}(A) = 1$. The author thanks
Prof N.V. Trung for this example.
\begin{example}
\label{2dimnotgcm}
Let $A = \mathbb{Q}[[s^4, s^3t, s^5t^3, t^4]]$. Using CoCoA \cite{cocoa} we get 

\noindent $A \cong \mathbb{Q}[[x,y,z,w]]/\q$ where $\q = (-x^2w + yz , -y^3+ xz, xy^2w - z^2)$.
So $G_{m}(A) \cong \mathbb{Q}[X,Y,Z,W]/\q^*$ where 
$\q^* = (-Z^2,YZ,XZ, -Y^4 + X^3W)$.  Set $R =\mathbb{Q}[X,Y,Z,W]$.   Using Singular \cite{singular} we
get $\dim \Ext^{3}_R(G_{\m}(A),R) = 1$. From graded local duality \cite[3.6.19]{BH} it follows that 
$H^{1}(G_{\m}(A))$ \emph{does not} have finite length. 
\end{example}

\begin{remark}
We do not have an example of a \CM \ local ring of dimension 2
such that $G_{\m}(A)$ is \GCM \ but $\xi_{\m}(A) =1$.
\end{remark}
%7777777777777777777
%section Observation
\section{An Observation}
\s In this section the setup is as given in \ref{setup}.

Let $x $ be $M$-superficial \wrt \  $I$. There exists  $c_i \in A$ such that $x
= \sum_{i= 1}^{m}c_ix_i$. Set $X = \sum_{i= 1}^{m}c_iX_i \in S = A[X_1,\ldots,X_m]$ and  $N
= M/xM$. We have the following commutative diagram:
\[
\xymatrixrowsep{3pc}
\xymatrixcolsep{3pc}
\xymatrix{
L^I(M)
\ar@{->}[r]^{\Pi}
\ar@{->}[dr]_{\Psi_x}
&L^I(M)(-1)
\ar@{->}[d]^{\Psi_X }
\\
&L^I(M)
}
\]
Notice $\Psi_x$ is multiplication by $x$ (we  consider $x$ as a degree zero element in $S$.)
\s
\label{cdiahom}
 Since $H^i(-)$ is a functor we have
the following commutative diagram:
\[
\xymatrixrowsep{3pc}
\xymatrixcolsep{3pc}
\xymatrix{
H^i\left(L^I(M) \right)
\ar@{->}[r]^{H^i(\Pi)}
\ar@{->}[dr]_{H^i(\Psi_x) }
&H^i\left(L^I(M)(-1)\right)
\ar@{->}[d]^{H^i(\Psi_X)}
\\
&H^i\left(L^I(M)\right)
}
\]
Notice that $H^i(\Psi_x) \colon H^i\left(L^I(M) \right) \rt H^i\left(L^I(M) \right)$ is just
multiplication by $x$ in each degree.
The following observation is central to all further results.
\begin{observation}
\label{cuteobs}
Suppose for some $n$ and some $i$ with $0 \leq i \leq \depth M -1$ the maps
\begin{align*}
H^i(\Pi)_n &\colon H^i\left(L^I(M) \right)_n \xar H^i\left(L^I(M) \right)_{n-1} \ \ \text{and} \\
H^i(\Psi_X)_n &\colon H^i\left(L^I(M) \right)_{n-1} \xar H^i\left(L^I(M) \right)_{n}  \ \ \text{are injective},
\end{align*}
 then $H^i(L^I(M))_n = 0$.
\end{observation}
\begin{proof}
Using \ref{cdiahom} we get that$ H^i(\Psi_x)_n \colon H^i\left(L^I(M) \right)_{n} \rt H^i\left(L^I(M) \right)_{n}$ is injective. Since $\lambda \left(H^i(L^I(M))_{n}\right)$ is finite we get that
$ H^i(\Psi_x)_n$ is an isomorphism.
 Therefore

\noindent $ H^i\left(L^I(M) \right)_{n} = x H^i\left(L^I(M) \right)_{n}$.
So $H^i(L^I(M))_n = 0$ by  Nakayama Lemma  .
\end{proof}

%888888888888888888888888888888

\section{A criteria for $\xi_I(M) \geq 2$}
\s In this section the setup is as given in \ref{setup}.

\s
\label{suffxi2}
Using Theorem \ref{asyym} and Proposition \ref{0lc}.2  we get $\xi_I(M) \geq 2$ if and only if
\[
H^1\left(L^I(M)\right)_{-1} = 0 \quad \text{and} \quad \lambda\left(H^1(L^I(M))\right) < \infty.
\]
By using the first of these conditions we were able to construct examples with $\xi_I(M) = 1$.
However we need a slightly different criteria if we want to show $\xi_I(M) \geq 2$. We do
this in the following:
\begin{proposition}
\label{xi2firstprop}
(with hypothesis as in \ref{setup}.)
 The following are equivalent:
\begin{enumerate}[\rm 1.]
\item
$\xi_I(M) \geq 2$.
\item
$H^1\left(L^I(M)\right)_{-1} = 0$.
\item
$H^1\left(L^I(M)\right)_n = 0$ for all $n <0$.
\item
$H^1(G_I(M))_{-1} = 0$.
\end{enumerate}
\end{proposition}
\begin{proof}
Using \ref{AtoA'} and \ref{basechangeC} we may assume that the residue field of $A$ is infinite.
Let $x $ be $M$-superficial \wrt \  $I$. There exists  $c_i \in A$ such that $x
= \sum_{i= 1}^{m}c_ix_i$. Set $X = \sum_{i= 1}^{m}c_iX_i \in A[X_1,\ldots,X_m]$ and  $N
= M/xM$.

$1. \implies 2.$ This is clear from \ref{suffxi2}.

$2. \implies 3.$ To see this we use part of (\ref{longH}),
\begin{equation*}
H^0\left(L^I(N)\right)_n  \xar H^1\left(L^I(M)\right)_{n-1}  \xrightarrow{H^1(\Psi_X)_n} H^1\left(L^I(M)\right)_n  \xar. \tag{*}
\end{equation*}
Since $H^0\left(L^I(N)\right)_n  = 0$ for $n < 0$ we get that
\[
\text{if} \  H^1\left(L^I(M)\right)_{-1} =0 \quad  \text{then} \quad
 H^1\left(L^I(M)\right)_n = 0 \ \text{for all} \ \  n < 0.
\]

$3. \implies 1.$  This follows from \ref{suffxi2} and the fact that
 $\lambda\left(H^1(L^I(M)\right)_n$ is finite for all $n$ and $ H^1\left(L^I(M)\right)_n = 0$ for
all $n \gg 0$.
Thus the conditions 1. 2. and 3. are equivalent.

$2. \implies 4.$ Set $L = L^I(M)$. We use part of the long exact sequence corresponding to (\ref{dag}),
\begin{equation*}
H^0\left(L\right)_{n-1}  \rt H^1\left(G_I(M)\right)_{n}  \rt H^1\left(L\right)_{n}
 \xrightarrow{H^1(\Pi)_n}  H^1\left(L\right)_{n-1}.  \tag{**}
\end{equation*}
Since $H^0\left(L\right)_{-2} = 0$ we get the result.

$4. \implies 2.$ If $ H^1\left(G_I(M)\right)_{-1} = 0$ then using (**) we get $H^1(\Pi)_{-1}$ is injective.
 Since  $H^0\left(L^I(N)\right)_{-1} = 0$ we get by (*) that  $H^1(\Psi_X)_{-1}$ is injective.
So by \ref{cuteobs} we get
$H^1\left(L^I(M)\right)_{-1} = 0$.
\end{proof}

In general one cannot compute $H^1(G_I(M))$. However we have the following noteworthy
special case:
\begin{theorem}\label{geq2b}
Let $(A,\m)$ be a Noetherian local ring with depth $\geq 2$ and residue field $k$. Let
 $G_{\m}(A) = R/\q$, where $R = k[X_1,\ldots,X_s]$. Then
\[
\Ext^{s-1}_{R}(G_{\m}(A),R)_{-(s-1)} = 0 \quad \ff \quad \depth G_{\m^n}(A) \geq 2 \ \text{for all} \
n \gg 0.
\]
\end{theorem}
\begin{proof}
Set $G = G_{\m}(A)$.
Notice $R$ is $*$-complete and  $R(-s)$ is the $*$-canonical module of $R$ cf. \cite[3.6.15]{BH}. Using the
local duality theorem for graded modules \cite[3.6.19]{BH} we get a homogeneous isomorphism
\[
H^1(G)^{\vee} \cong \Ext^{s-1}_R\left(G,R(-s)\right).
\]
For any graded $R$-module $E$ we have $(E^{\vee})_i = \Hom_{k}(E_{-i},k)$ \cite[Section 3.6]{BH}.
Notice
\begin{align*}
\lambda\left(H^1(G)_{-1}\right) &= \lambda\left( H^1(G)^{\vee}_1 \right) \\
 &=  \lambda\left(\Ext^{s-1}_R\left(G,R(-s)\right)_{1}\right) \\
                           &=  \lambda\left(\Ext^{s-1}_R\left(G,R\right)_{-(s-1)}\right).
\end{align*}
The result follows by using Proposition \ref{xi2firstprop}.
\end{proof}

\section*{Acknowledgements}
Its a pleasure to thank Prof M. Brodmann, Prof J.K. Verma, Prof
M.E. Rossi, Prof N.V. Trung and Dr. A.V. Jayanthan for many discussions.

\providecommand{\bysame}{\leavevmode\hbox to3em{\hrulefill}\thinspace}
\providecommand{\MR}{\relax\ifhmode\unskip\space\fi MR }
% \MRhref is called by the amsart/book/proc definition of \MR.
\providecommand{\MRhref}[2]{%
  \href{http://www.ams.org/mathscinet-getitem?mr=#1}{#2}
}
\providecommand{\href}[2]{#2}

\end{document}